\input amssym.def
\input amssym.tex%
\magnification = 1200

\pageno=1
\mag=1200 
\baselineskip=5mm
 \hoffset=-4mm
  \voffset=8mm
 \hsize=14cm
\vsize=190mm
\input xy
\xyoption{all}
\par
\topskip = 10 truemm

\overfullrule=0pt

\rightskip=10mm

\def \noi {\noindent}
\def \dis {\displaystyle}

\def \cqfd{\unskip\kern 6pt\penalty 500
 \raise -2pt\hbox{\vrule\vbox to10pt{\hrule width 4pt
\vfill\hrule}\vrule}\par}                

\def\adots{\mathinner{\mkern2mu\raise1pt\hbox{.}
\mkern3mu\raise4pt\hbox{.}\mkern1mu\raise7pt\hbox{.}}}
\def\pmb#1{\setbox0=\hbox{#1}%
\kern-.025em\copy0\kern-\wd0
\kern.05em\copy0\kern-\wd0
\kern-.025em\raise.0433em\box0}

\font \titre = cmbx10 scaled 1400
\font \tite=cmbxti10 scaled 1400

\font \bfdouze = cmbx12


\def\build#1_#2^#3{\mathrel{\mathop{\kern 0pt#1}\limits_{#2}^{#3}}} 

\def\diagram#1{\def\normalbaselines{\baselineskip=0pt
\lineskip=8pt\lineskiplimit=1pt} \matrix{#1}}
\def\build#1_#2^#3{\mathrel{\mathop{\kern 0pt#1}\limits_{#2}^{#3}}} 
\def\hfl#1#2{\smash{\mathop{\hbox to 18mm{\rightarrowfill}}\limits^{\scriptstyle#1}_{
\scriptstyle#2}}}
\def\vfl#1#2{\llap{$\scriptstyle #1$}\left\downarrow \vbox to
4mm{}\right.\rlap{$\scriptstyle #2$}}

\def\boxit#1#2{\setbox1=\hbox{\kern#1{#2}\kern#1}%
\dimen1=\ht1 \advance\dimen1 by #1 \dimen2=\dp1 \advance\dimen2 by #1
\setbox1=\hbox{\vrule height\dimen1 depth\dimen2\box1\vrule}%
\setbox1=\vbox{\hrule\box1\hrule}%
\advance\dimen1 by .4pt \ht1=\dimen1
\advance\dimen2 by .4pt \dp1=\dimen2 \box1\relax}

\input xy

\headline={\ifnum\pageno=1 {\hfill} \else{\hss \tenrm -- \folio\ -- \hss} \fi}
\footline={\hfil}



\centerline{\titre INCIDENCE DIVISOR }
\vskip 1cm

\centerline{\bfdouze D. Barlet\footnote{\rm{(*)}}{\rm{Senior chair of Complex Analysis and Geometry in Institut Universitaire de France}} and \bfdouze{ M. Kaddar}}

\vskip 0,5cm

\centerline{\it Universit\'e Henri Poincar\'e Nancy 1}

\centerline{\it Institut Elie Cartan}

\centerline{\it UMR 7502 CNRS - INRIA - UHP} 

\centerline{\it BP 239 - F - 54506 Vand\oe uvre-l\`es-Nancy Cedex}

\vskip 1,5cm
\noi {\tite 0. Introduction.}

\smallskip

\noi The main purpose of this paper is to prove an important generalization of the construction of
the Incidence Divisor given in [BMg1] in the case of an ambient manifold.

\noi Let us first recall briefly the setting : let $Z$ be a complex manifold and $\;(X_s)_{s
\in S}\;$ an analytic family of (closed) $n-$cycles in $Z$ parametrized by a reduced complex
space $S$.

\noi To a $(n+1)-$codimensional subspace $Y$ in $Z$, which is assumed to be a locally  complete
intersection and  to satisfy the following condition :

\smallskip

\item{(C1)} $\;$the analytic set $\;(S \times |Y|) \cap |X|\;$ in $S \times Z$\footnote{(**)}{$X$ denotes
the graph of the family $(X_s)_{S \in s}$, which is a cycle in $S \times Z$ and $|X|$ the
support of this cycle.} is $S$ -proper and finite on its image $|\Sigma_Y|$ which is nowhere
dense in $S$,

\smallskip

\noi an effective Cartier divisor $\Sigma_Y$ in $S$, called the "incidence divisor of $Y$ in
$S$", is defined with support $|\Sigma_Y|$, and nice functorial properties of this construction
are proven. Of course, no assumption is made on the  singularities of $S$.

\noi A relative version is also given : when $Y$ varies in a flat family over $T$ in
such a way that (C1) remains  true, $\Sigma_Y$ moves in a flat family over $T$. As a
consequence, $\Sigma_Y$ depends only on the underlying cycle of $Y$ for a connected flat
deformation of a locally complete intersection (nilpotent) structure inducing a fixed cycle.

\noi But the general invariance question was not solved in [BMg1].

\smallskip

\item{[\bf{Q1}]}\rm Does  the Cartier divisor $\Sigma_Y$ only depend on the cycle underlying the 
locally complete intersection ideal of ${\cal O}_Z$ defining $Y$ ?

\noi Of course another natural (but stronger) question arises if [\bf{Q1}]\rm admit a positive
answer:

\smallskip

\item{[\bf{Q2}]}\rm Is it possible to construct the incidence divisor $\Sigma_Y$ (Cartier and
effective) for any cycle $Y$ satisfying (C1) (with nice functorial properties) ?

\smallskip

\noi Of course the construction of a line bundle given in [K] suggested that the cohomological
method could help to solve these questions. Also the hypothesis used in [K] was weaker
than (C1). In fact the condition (C1) has a nice geometrical meaning and is very 
common {\bf generically} on $S$. But in concrete examples some degenerations appear very
often. 

\smallskip

\noi The following weaker condition seems much more attractive for applications :

\smallskip

\item{(C2)} the analytic set $\;\big(S \times |Y|\big) \cap |X|\;$ in $S \times Z$ is 
$\;S-$proper and {\bf generically} finite on its image $|\Sigma_Y|$ which is nowhere 
dense in $\;S$.

\smallskip

\noi But extending a Cartier divisor through an analytic set is, in general hopeless 
without strong assumptions on the singularities of $S$. So the geometric method 
developed in [BMg1] cannot treat this (C2) situation in general.

\smallskip

\noi Our generalization gives the following positive answer :

\smallskip

\noi {\bf [Q2] has a positive answer with the assumption (C2) (see the Theorem  in \S  2).

\noi\rm In fact, using the cohomological tools built in [K], that is to say the existence of a relative fundamental class in Deligne-Beilinson cohomology for the analytic family
$(X_s)_{s \in S}$, $\;$ we are able to prove also a relative version for [Q2], which solves also
[Q1], assuming only (C2).}

\smallskip

\noi To conclude this article, we prove a formula computing the intersection number of the
incidence divisor $\Sigma_Y$ with a curve in $S$ in terms of intersection numbers of $Y$ in
$Z$. This is useful, for instance, to compute the multiplicities of the incidence divisor,
but can be used also in a more global setting ($Z$ projective) to obtain a simple criterion
of nefness for the incidence divisor. Combined with results in [BMg2] this may lead to  vanishing theorems on $S$.

\noi Finally, let us remark that this article is a modest contribution, in the local analytic
context and at the geometric level, to P. Deligne's  program \big(see [D]\big), which has been
developed in a global projective algebraic setting (with rather strong hypothesis) by R.
Elkik in [E].

\smallskip

\noi In \S 1 we explain how to build up a Cartier divisor from "special" cohomology classes in $I\!\!H_T^2 (S, {\Bbb Z} (1)_{\cal D})$
following in our local analytic setting some simple
ideas from algebraic geometry \big(more precisely pseudo-divisors as described in 
[F]\big). Of course, essential singularities have to be avoided.

\smallskip

\noi In \S 2 we use the cohomological tools from [K] to build  a cohomology class in $I\!\!H^2_T \big(S, {\Bbb Z}(1)_{\cal D}\big) \;$ and we prove that it  satisfies the criterion
proved in \S 1.

\smallskip

\noi In \S 3 we prove various functorial properties of our construction. This leads to the
comparaison with [BMg 1] and give effectivity of our incidence divisor. The relative case
and the positive answer to [Q1] conclude \S 3.

\smallskip

\noi In \S  4 we compute the intersection number of the incidence divisor with a curve in
$S$.

\bigskip

\noi Finally, we have postponed till the Appendix some basic properties of (relative) 
fundamental classes in Deligne cohomology such as compatibility with intersection, projection
formula and traces, for which it seems that no reference is avaible in our local analytic
context.

\vskip 1cm

\noi {\tite 1. A  criterion.}

\smallskip

\noi Let $S$ be a reduced complex analytic space and $\;T \subset S\;$ a closed analytic
subset of $S$ with empty interior. Denote by $\; j:S-T\hookrightarrow S\;$ the inclusion
and by $\;{\cal O}^*_S\;$ the sheaf of abelian (multiplicative) groups of non vanishing
germs of holomorphic functions on $S$. We shall consider  $\;{\cal O}^*_S \;$ as a subsheaf
of the sheaf $\;m^*_T\;$ of germs of meromorphic functions (not vanishing  identically  on any open set
in $S$) with zeros and poles contained in $T$. So $\; m^*_T |_{S-T} = {\cal O}^*_S|_{S-T}$.

\smallskip

\noi There is a canonical isomorphism of sheaves 
$$j_* j^* {\cal O}^*_S \big/{\cal O}^*_S \build \longrightarrow_{}^{\sim} \underline{H}^1_T
({\cal O}^*_S)$$

\noi and a natural inclusion
$$m^*_T \hookrightarrow j_* j^* {\cal O}^*_S.$$

\noi So the quotient sheaf $\; m^*_T\big/{\cal O}^*_S \;$ is naturally embedded as a subsheaf
of $\; \underline{H}^1_T ({\cal O}^*_S)$. We shall denote it by $\; \underline{H}^1_{[T]}
({\cal O}^*_{S})$.

\smallskip

\noi Now let us consider the sheaf $\; \underline{C}_T \;$ of Cartier divisors in $S$ with
support in $T$.

\noi This is a sheaf of abelian groups on $S$ supported in $T$ and we have a natural homomorphism
of sheaves of abelian groups
$$ {\rm div} : m^*_T \longrightarrow \underline{C}_T$$

\noi defined  by ${\rm div} \dis \Big({f \over g}\Big) = \{f=0\}-\{g=0\} \;$ for $\;f,g \in
{\cal O} (U)\;$ such that $\; \{f=0\}$ and $\{g=0\}\;$ are contained in $\;U \cap T$.

\noi The kernel of div is clearly $\; {\cal O}^*_S \subset m^*_T \;$ and so we obtain an
isomorphism of sheaves of abelian groups
$${\rm Div} : \quad \underline{H}^1_{[T]}  \; ({\cal O}^*_S) \build \longrightarrow_{}^{\sim} 
\underline{C}_T.$$

\noi The reader may compare with [F] p. 31 for the algebraic version of these facts. We
want to give a simple criterion in order that a section $\; \sigma \;$ of $\; \underline{H}^1_T
({\cal O}^*_S)\;$ be  a section of the subsheaf $\; \underline{H}^1_{[T]} ({\cal O}^*_S)$, and so, produces a Cartier divisor via Div.

\smallskip

\noi Consider the normalisation map $\nu : \widetilde{S} \longrightarrow S \;$ for $S$ and
the following sheaf homomorphisms  
$$\matrix{
d \; {\rm log} & : {\cal O}^*_S \longrightarrow \Omega^1_S \hfill & (g \longrightarrow 
dg/g) \cr
\noalign{\vskip 0,3cm}
\nu^* & : \Omega^1_S \longrightarrow \nu_* \Omega^1_{\widetilde{S}} \hfill & \cr
\noalign{\vskip 0,3cm}
q &: \nu_* \Omega^1_{\widetilde{S}} \longrightarrow \nu_* \omega^1_{\widetilde{S}} \hfill
& \cr}$$

\smallskip

\noi where $\; \omega^1_{\widetilde{S}} \;$ is the sheaf defined in [B.3].

\smallskip

\noi The normality of  $ \widetilde{S} $ gives the isomorphism
$$\omega^1_{\widetilde{S}}\simeq i_* i^* \Omega^1_{\widetilde{S}} $$

\noi where $\; i : \widetilde{S} - {\rm Sing}(\widetilde{S}) \hookrightarrow \widetilde{S} \;$ 
is  the inclusion of the smooth points in $\; \widetilde{S}$.

\smallskip

\noi Let $\widetilde{T} := \nu^{-1} (T). \;$ Consider now the sheaf morphism
$$\varphi : \underline{H}^1_T ({\cal O}^*_S) \longrightarrow \underline{H}^1_T (\nu_* 
\omega^1_S) \simeq \nu_* \underline{H}^{1}_{\widetilde{T}} (\omega^1_S) $$

\noi deduced from $d$Log,  $\nu^*$ and $\;q$.

\smallskip

\noi Denote by $\; \underline{H}^1_{[1.\widetilde{T}]} \big(\omega^1_{\widetilde{S}}\big)
\;$ the subsheaf of the sheaf $\; \underline{H}^1_{\widetilde{T}} \big(\omega^1_{\widetilde{S}}
\big) \;$ of sections annihilated by the reduced ideal $\; {\cal J}_{\widetilde{T}}\;$ of
$\; \widetilde{T}\;$ in $\;\widetilde{S}$. Note that it is equivalent for a section 
$\;\sigma\;$ of $\; \underline{H}^1_{\widetilde{T}} \big(\omega^1_{\widetilde{S}} \big)
\;$ to be annihilated by $\;{\cal J}_{\widetilde{T}} \;$ or to be annihilated by $\;
{\cal J}_{\widetilde{T}} \;$ at generic points of codimension 1 of $\;\widetilde{T}\;$
because the sheaf $\; \omega^1_{\widetilde{S}} \;$ has no torsion and satisfies the analytic
continuation property in codimension 2 \big(see [B.3]\big). But this latter condition is
very simple to check because at generic points of codimension one of $\;\widetilde{T},\;
\widetilde{S}\;$ is smooth and $\; \widetilde{T}\;$ is a smooth hypersurface. So we are
checking a pole of order one for a singularity of a  holomorphic 1-form along a smooth
hypersurface.

\smallskip

\noi Our criterion is given by the following

\smallskip

\noi {\bf Proposition 1.}

\noi \it{In the above notation, we have}
$$\varphi^{-1} \big[\nu_* \underline{H}^1_{[1.\widetilde{T}]} \big(\omega^1_{\widetilde{S}}
\big)\big] = \underline{H}^{1}_{[T]} ({\cal O}^*_S)$$\rm
\vfill\eject
\noi {\bf Proof.}

\noi First, if $\;g \in m^*_T \;$ we have $\; \varphi [g] = \nu_* \Big(\nu^*\Big({dg \over
g}\Big)\Big)\;$ in $\; \nu_* \underline{H}^1_{\widetilde{S}} (\omega^1_{\widetilde{S}}) 
\;$ and we want to see that $\;\nu^* \Big({dg \over g} \Big)\;$ is annihilated at generic
codimension 1 points of $\; \widetilde{T}\;$ in $\; \widetilde{S}\;$ by ${\cal J}_{\widetilde{T}}
.\;$ This is clear. Conversely if $\; g \in j_* j^* {\cal O}^*_S\;$ we want to show that
if $\; \nu^* \Big({dg \over g}\Big)\;$ has at most simple poles at generic codimension 1
points in $\; \widetilde{T} \;$ then $\;g \in m^*_T.\;$ This is a consequence of the following
classical Lemma

\smallskip

\noi {\bf Lemma 1.}

\noi {\it Let $\;D^*:=\{z \in {\Bbb C} / 0<|z|<1\}\;$ and let $\; g \in {\cal O}(D^*).\;$ 
Assume that $\; z \dis {g' \over g} \in {\cal O} (D). \;$ Then $\;g\;$ is meromorphic at
$0$.}
\rm
\smallskip

\noi {\bf Proof.}

\noi Denote by $\;h\;$ the holomorphic function $\dis z {g'\over g}\;$ on $\;D.\;$ Then,
on $\;D^*,g\;$ satisfies the differential equation $\; g'= \dis {h \over z} \cdot g.\;$
This equation has a simple pole at $\;z=0$. So $\;g\;$ has moderate growth at 0. The
conclusion follows. $\hfill \blacksquare$

\vskip 1cm

\noi Using this Lemma we obtain that $\; \nu^*g\;$ is meromorphic on $\; \widetilde{S}\;$
minus an analytic set of codimension $\ge 2\;$. So $\; \nu^*g \;$ is meromorphic on $\;\widetilde{S}\;$
and also on $\;S$. $\hfill \blacksquare$

\vskip 1cm

\noi {\tite 2. The main construction.}

\smallskip

\noi Now we shall consider the following situation : 

\noi Let $\;Z\;$ be a complex manifold of dimension $\;n+p\;$ at let $\;(X_s)_{s \in S}\;$
be an analytic family of closed $\;n-$cycles in $\;Z\;$ parametrized by a reduced complex
space $\;S.\;$ Let $\;Y\;$ be a cycle of pure dimension $\;(p-1)\;$ in $\;Z\;$ and assume
the following condition :
$$ [AP]\footnote{(*)}{$AP\;$ for admissible pole} \left[ \matrix{
& {\rm Let} \; X \subset S\times Z \; \hbox{be the graph of the family} \; (X_s)_{s \in S} \;\hbox{,}
 \; p_X : X \longrightarrow Z \hfill \cr
& \hbox{ the projection induced by the canonical projection of }\;S\times Z\;\hbox{on} \;Z\;\hbox{ and let} \hfill \cr
& p_X^{-1}(|Y|):=(S\times |Y|) \cap |X|. \; \hbox{We require  that :} \hfill\cr
\noalign{\vskip 0,2cm}
& 1) \; \hbox{codim} \; p_X^{-1} (|Y|) \; {\rm in} \; S \times Z \;{\rm is}\;\; n+p+1 \; \hbox{(the 
expected codimension)} \hfill\cr
\noalign{\vskip 0,2cm}
& 2) \hbox{the restriction of} \;\; \pi : |X| \longrightarrow S \;\;{\rm to}\;\; p_X^{-1}
(|Y|) \; {\rm is} \hbox{ proper and generically } \hfill\cr
& \hbox{finite on its image} \quad T:=\pi(p_X^{-1}|Y|) \hfill }\right. $$

\noi {\bf Remark} $\quad$ Under $\;[AP]\;$ one can define the intersection in $\; S \times
Z\;$ of the cycles $\; p_{X}^* Y := X \build ._{p}^{} Y \;$ where $\; p:S\times Z \longrightarrow Z\;$ is the
second projection\footnote{(**)}{See [B.1] ch 6 \S 3} and also the direct image $\;\pi_*
(p_X^*Y)\;$ as a cycle. We obtain in this way a codimension 1 Weil cycle in $\;S\;$ with
support in $\;T.\;$ Everything comes down to showing  that this "rough" geometric construction
underlies the construction of a natural Cartier divisor in $\;S\;$ generalizing the results
of [BMg1].

\smallskip

\noi We shall prove the following result :

\smallskip

\noi {\bf Theorem .}

\noi {\it In the situation described above to each cycle $\;Y\;$ satisfying $\;[AP]\;$ we
associate a Cartier divisor $\; \Sigma_Y\;$ in $\;S\;$ with support in $\; T:= \{ s \in S/
|Y| \cap |X_s| \neq \emptyset \}\;$.{\it This construction satisfies }
\smallskip

\noi {\it $1$) For $\; Z={\Bbb C}, \; S={\rm Sym}^k({\Bbb C}) \simeq {\Bbb C}^k/\sigma_k, \;\;
(X_s)_{s \in S} \;$ the universal family of 0-cycles of degree $\;k\;$ in $\;{\Bbb C}\;$
parametrized by $\;S\;$ and $\;Y=1.\{0\}\;$, we have $\;\Sigma_Y =\{s_k =0\}\;$ where $\;s_k
: {\rm Sym}^k({\Bbb C}) \longrightarrow {\Bbb C} \;$ is the $k-$th elementary symetric
function.

\smallskip\noindent The map $\; Y\rightarrow {\Sigma_{Y}}\;$  is compatible with\smallskip

\noi $2$)allowed base change (so is local in $\;S$)

\smallskip

\noi $3$) localization in $\;Z$

\smallskip

\noi $4$) direct image of cycles

\smallskip

\noi $5$) inverse image of cycles

\smallskip

\noi $6$) analytic deformation of $\;Y$ as a cycle.

\smallskip

\noi Moreover these properties characterize our construction.}\rm
\bigskip\noindent Precise formulations of conditions
$2$) to $6$) will be given at the beginning of \S 3. 

\noi We shall now give the definition of $\;\Sigma_Y\;$ using tools from [K] and the criterion
of the previous paragraph. The functorial properties 2) to 6) will be proved in \S 3.

\smallskip

\noi Let $\;C^{S\times Z}_{X/S}\;$ be the relative fundamental class in Deligne cohomology
for the family $\;(X_s)_{s \in S}.\;$ Recall that $\;C^{S\times Z}_{X/S} \;$ is a (global)
section of the sheaf
$$\underline{I\!\!H}^{2p}_{|X|} \big({\Bbb Z} (p)_{{\cal D}/S} \big)$$

\noi where we denote by $ \Bbb Z(p)_{{\cal D}/S}$ the complex of sheaves on $ S \times
Z$
$${\Bbb Z}(p)_{{\cal D}/S} : = (2i\pi)^{p}( {\Bbb Z}) \longrightarrow {\cal O}_{S \times Z/S}{\longrightarrow}{\Omega}^{1}_{S \times Z/S} \longrightarrow\cdots {\longrightarrow}{ \Omega}^{p-1}_{S \times
Z/S} $$
\noindent with the differential $S$- relative (for the construction see [K] Th. II).

\smallskip

\noi Now let us consider $(Y_s)_{s\in S}$  the constant family of cycles in $Z$ given  by  $ Y_s=Y$ for all $s$ in  $S$. We have also a relative fundamental class $C^{S \times Z}_{S \times Y/S} $
which is the section of the sheaf $ \underline{I\!\!H}^{2(n+1)}_{S \times |Y|}  \big({\Bbb Z}(n+1)_{{\cal D}/S})
\big)$  obtained from $C^Z_Y$ \big(a section of $ \underline{I\!\!H}^{2(n+1)}_{|Y|}
({\Bbb Z}(n+1)_{{\cal D}})$ on $Z$\big) by pull back and the natural map
$${\Bbb Z}(n+1)_{\cal D} \longrightarrow {\Bbb Z}(n+1)_{{\cal D}/S}$$
\noi deduced from the quotients $ \Omega^{\bullet}_{S \times Z} \longrightarrow \Omega^{\bullet}_{S \times
Z/S}$.
\noi Now the cup-product $ C^{S \times Z}_{X/S} \cup C^{S \times Z}_{S \times Y/S} $ is
a section of the sheaf
$$\underline{I\!\!H}^{2(n+p+1)}_{p_{X}^{-1}(|Y|)} \big({\Bbb Z}(n+p+1)_{{\cal D}/S}\big)$$ ( for the definition of the cup product in the absolute case in Deligne cohomology see [E.V]; the relative case is similar).
 Using the trace map \big(see [K] \S V\big)
$$Tr_{/S} : \pi_* \underline{I\!\!H}^{2(n+p+1)}_{p_{X}^{-1}(|Y|)} \big({\Bbb Z}(n+p+1)_{{\cal D}/S}\big)
\longrightarrow \underline{I\!\!H}^2_T \big({\Bbb Z}(1)_{\cal D}\big)$$

\noi we obtain a section $ \alpha$ of the sheaf $ \underline{I\!\!H}^2_{T} \big({\Bbb Z}(1)_{\cal D}\big)$ on $S.$ But we have 
$${\Bbb Z}(1)_{\cal D} : = 2i\pi {\Bbb Z} \longrightarrow {\cal O}_S $$

\noi and so $ {\Bbb Z}(1)_{\cal D}$ is quasi-isomorphic to the sheaf ${\cal O}^*_S[-1]$
via the exponential. We shall then consider $ \alpha$ as a section of the sheaf
$ \underline{H}^1_T ({\cal O}^*_S)$ and our goal is now to prove that we have in fact a
section of the subsheaf $ \underline{H}^1_{[T]} ({\cal O}^*_S) $ which is isomorphic
via Div (see \S 1) to the sheaf $ \underline{C}_T$ of Cartier divisors supported in
$T$.

\noi So the key to our construction is to show that, near the generic points in 
$\widetilde{T}$ of codimension 1 in $\widetilde{S}$ (they are all of codimension one by
our previous remark !), the image of $\alpha$ by $\varphi$  has at
most a simple pole along $\widetilde{T}$.

\smallskip

\noi Then using the compatibility of our construction with base change as long as $Y$
satisfies $[AP],$ we are reduced to show that, for  $S$ {\bf smooth},  $\alpha$ is given by our "rough" geometric construction which could be seen a follows :

\smallskip

\noi Let $ C^{S \times Z}_X$ be the absolute fundamental class of the cycle $X$ in
$S \times Z,$ seen  as a section of the sheaf $ \underline{I\!\!H}^{2p}_{|X|} \big(
{\Bbb Z}(p)_{\cal D}\big)$.

\noi Let $ C^{S \times Z}_{S \times Y} $ be the absolute fundamental class of the 
cycle $ S \times Y$ in $ S \times Z$.

\smallskip

\noi Now the intersection $ p_X^* Y : = X . (S \times Y) $ is well defined as a cycle
in $ S \times Z$ 

\noi($S$ is smooth) and the absolute fundamental class of $p^{*}_{X} Y$ is
the cup product 
$$ C^{S \times Z}_{X} \cup C^{S \times Z}_{S \times Y} $$

\noi as a section of the sheaf $ \underline{I\!\!H}^{2(n+p+1)}_{p_X^{-1}(|Y|)}  \big(
{\Bbb Z}(n+p+1)_{\cal D}\big)$ (see  Appendix II.A).

\smallskip

\noi Of course via the obvious map
$$\underline{I\!\!H}^{2(n+p+1)}_{p_X^{-1}(|Y|)} \big({\Bbb Z}(n+p+1)_{\cal D}\big) \longrightarrow
\underline{I\!\!H}^{2(n+p+1)}_{p_X^{-1}(|Y|)} \big({\Bbb Z}(n+p+1)_{{\cal D}/S}\big)$$

\noi this cup product is sent to $ C^{S \times Z}_{X/S} \cup C^{S \times Z}_{S \times Y/S}$.

\noi  We have a commutative diagram of sheaves on $S$:
$$\diagram{
\pi_* \underline{I\!\!H}^{2(n+p+1)}_{p_X^{-1}(|Y|)} \big({\Bbb Z}(n+p+1)_{\cal D}\big) &
\hfl{\lambda}{} & \pi_* \underline{H}^{n+p+1}_{p_X^{-1}(|Y|)} (\Omega^{n+p+1}_{S \times Z})\cr
\vfl{}{Tr} && \vfl{}{\pi_*} \cr
\underline{I\!\!H}^2_T \big({\Bbb Z}(1)_{\cal D}\big) & \hfl{{1 \over 2i\pi} d {\rm log}}{} &
\underline{ H}^1_T (\Omega^1_S) \cr } \leqno (D)$$

\noi where the vertical arrows are absolute traces induced by integration maps 
which give, for $\pi_{*}$, values in the sheaf of the holomorphic forms on $S$ assumed to be smooth  \big(see [K]\big) and where $\lambda$
is deduced from the boundary map
$$\partial:  \underline{I\!\!H}^{2k}_A \big({\Bbb Z} (k)_{\cal D}\big) \longrightarrow \underline{I\!\!H}^{2k+1}_A
\big(\Omega^k [-k-1]\big)  $$

\noi of the exact sequence of complexes 
$$0 \longrightarrow \Omega^k [-k-1] \longrightarrow {\Bbb Z}(k+1)_{\cal D} \longrightarrow
{\Bbb Z}(k)_{\cal D} \longrightarrow 0.$$

\noi For $k=1$ and $A=T$ in $S$, it is a simple exercice to see that  $\partial$
 gives $d$ Log\par\noindent \Big(up to $\dis {1 \over 2i\pi}$ !\Big) using the quasi-isomorphism
$${\Bbb Z}(1)_{\cal D} \simeq {\cal O}^*_S [-1]$$
and the exact sequence
$$0\rightarrow \Omega^{1}[-2]\rightarrow{\Bbb Z}(2)_{\cal D}\rightarrow {\Bbb Z}(1)_{\cal D}\rightarrow 0$$ 
\noi In fact the commutativity of $(D)$ reflects the compatibility between direct images of
cycles and the trace of their fundamental classes (in Deligne cohomology and in holomorphic
cohomology).

\smallskip

\noi Denote by $\widehat{T}$ the cycle $\pi_* (p_X^*Y).$ We have 
$${1 \over 2i\pi}d {\rm Log}  Tr(C^{S \times Z}_{p_X^*Y})= C^S_{\widehat{T}} \in 
\underline{H}^1_T (\Omega^1_S)$$

\noi and the fundamental class $C^S_{\widehat{T}}$ is annihilated by the reduced ideal
of $|\widehat{T}|$ in $S$ (still assumed smooth).

\smallskip

\noi So, to conclude that $\alpha$ is a section of $ \underline{H}^1_{[T]} ({\cal O}^*_S)$
it is enough to show that the following triangle commutes

$$\xymatrix{\pi_* \underline{I\!\!H}^{2(n+p+1)}_{p^{-1}_X(|Y|)} \big({\Bbb Z}(n+p+1)_{\cal D}\big)\ar[rr]\ar[rd]^{Tr}&&\ar[ld]^{Tr_{/S}}\pi_* \underline{I\!\!H}^{2(n+p+1)}_{p^{-1}_X(|Y|)} \big({\Bbb Z}(n+p+1)_{{\cal D}/S}\big) \\ 
& \underline{I\!\!H}^2_T \big({\Bbb Z}(1)_{\cal D}\big)&}$$

\noi for $S$ smooth. But, by definition of the trace $Tr$, it factorises through
$Tr_{/S}$ \big(see [K] p.320\big).

\noi So  $\alpha$ is a section of $ \underline{H}^1_{[T]} ({\cal O}^*_S)$ and
define via Div : $\underline{H}^1_{[T]} ({\cal O}^*_S) \longrightarrow \underline{C}_T$
a Cartier divisor in $S$ with support in $T,$ denoted by $\Sigma_Y$.

\bigskip

\noi It is clear from our construction that the Cartier divisor $\Sigma_Y$ in $S$ 
satisfies, as a codimension one  cycle in $S$ (i.e as Weil divisor)
$$\Sigma_Y = \pi_* (p_{X}^* Y)$$

\noi where $p_X^*Y:=X \build ._{p}^{} Y$ in $ S \times Z$\footnote{(*)}{see [B.1] ch.6 \S 3 for the
definition. Here $p:S \times Z \longrightarrow Z$ is the projection and $Z$ is
smooth}.

\smallskip

\noi But, of course, the proof that the Weil divisor $ \pi_*(p_X^*Y)$ is Cartier does not seem
 to be clear (for $S$ general) without using Deligne cohomology.\smallskip\noindent
\bf{Remark.}\rm\smallskip\noindent
With condition (C2), it is no hard to see that, for all $i\in {\Bbb N}$, [K] give us the following commutative diagram
$$\xymatrix{
&I\!\!H^{2n+i}_{|Y|}(Z, {\Omega^{*}_{Z}/F^{n+1}\Omega^{*}_{Z}})\ar[r]
\ar[d]&
I\!\!H^{2n+i+1}_{|Y|}(Z,{\Bbb Z}(n+1)_{{\cal D}})\ar[r]
\ar[d]&I\!\!H^{2n+i+1}_{|Y|}(Z,{\Bbb Z}(n+1))\ar[d]&&\\
&\Gamma(S, {\cal H}^{i}_{T} ({\cal O}_{S}))\ar[r]
&\Gamma(S, {\cal H}^{i}_{T}({\cal O}^{*}_{S}))\ar[r]&
\Gamma(S, {\cal H}^{i+1}_{T}( {\Bbb Z}(1))),&&}$$
where $F^{\bullet}$ is the shifted {\it{filtration b\^ete}}\rm\quad on the De Rham complex on $Z$.\par\noindent
This  gives, in  particular, the map
$$\Psi^{X}_{Z,S}: I\!\!H^{0}(Z, I\!\!{\cal H}^{2n+2}_{|Y|}(Z,{\Bbb Z}(n+1)_{{\cal D}}))\rightarrow
 \Gamma(S, {\cal H}^{1}_{T} ( {\cal O}^{*}_{S}))$$
\vskip 1cm

\noi {\tite 3. Functorial properties.}

\smallskip

\noi We shall begin this section with precise statements and proofs of Properties 2) to 
6) of the theorem.

\smallskip\bigskip

\noi {\bf 2) Base change.}

\noi We consider a holomorphic map $ \lambda : S' \longrightarrow S $ of reduced complex
space and we assume that $Y$ satisfies $[AP]$ for the analytic family $(X_s)_{s \in S}$
and also for the family $(X_{\lambda(s')})_{s' \in S'}.$ This last condition reduces to
the following hypothesis :

\noi Let $ R \subset T$ be the analytic subset in $T=\big|\pi_*\big(p_X^*(Y)\big)\big|$
where $\pi : |p_X^* Y| \longrightarrow T $ has positive dimensional fibers. Then $R$
is a closed analytic subset with empty interior in $T$ from $[AP]$ \big(relative to the
family $(X_s)_{s \in S}\big)$. Then $Y$ will satisfy $[AP]$ for the family
$(X_{\lambda (s')})_{s' \in S'\quad}$ iff
$$(H) \left \{ \eqalign{
& \lambda^{-1} (R)  \hbox{\bf has no interior points in}\quad  \lambda^{-1} (T)\cr
& \hbox{\bf which has no interior point in}\quad  S'. \cr} \right. $$  

\noi Under these hypothesis, we want to show that
$$\vbox{\boxit{5pt}{$\lambda^* \big(\sum_Y^S\big)=\sum_Y^{S'}$}}$$

\noi Observe  that the pull back of the Cartier divisor is well defined under $(H)$

\smallskip

\noi {\bf Proof.}

\noi Let $p':S' \times Z \longrightarrow Z$ be the projection. We want to prove the
commutativity of the square 
$$\diagram{
\pi_* \underline{I\!\!H}^{2(n+p+1)}_{p_X^{-1}(|Y|)} \big({\Bbb Z}(n+p+1)_{{\cal D}/S}\big) &
\hfl{\lambda^*}{} & \pi'_* \underline{I\!\!H}^{2(n+p+1)}_{p'^{-1}_X(|Y|)} \big({\Bbb Z} 
(n+p+1)_{{\cal D}/S}\big)\cr
\vfl{}{T_{/S}} && \vfl{}{T_{/S'}} \cr
\underline{I\!\!H}^2_T \big({\Bbb Z}(1)_{\cal D}\big) & \hfl{\lambda^*}{} &
\underline{I\!\!H}^2_{T'} \big({\Bbb Z} (1)_{\cal D}\big) \cr } $$

\noi where $T' := \pi' \big[p'^{-1}_X (|Y|)\big]$.

\smallskip

\noi This is consequence of the stability by base change of the relative trace in Deligne
cohomology (see Appendix I and [K]). $\hfill \blacksquare$
\vfill\eject
\noi {\bf 3) The construction is local on Z.}

\noi Let $U$ an open set in $Z$ and let
$$S_U:=\{s \in S \big/ |X_s| \cap |Y| \subset U\}.$$

\noi Then $S_U$ is open in $S$ and we want to show that
$$\vbox{\boxit{5pt}{$\sum_Y |_{S_{U}}=\sum_{Y|_U}$.}}$$

\noi Notice  that when $Y$ satisfies $[AP],  Y \cap U$ satisfies $[AP]$ for
the family $(X_s \cap U)_{s \in S_{U}}$.

\smallskip

\noi The proof is a consequence of the compatibility of our constructions with restrictions
to open sets. We have already put some emphasis on that by arguing at the level of sheaves
! $\hfill \blacksquare$

\vskip 1cm

\noi {\bf 4) Direct image of cycles.}

\noi Let $ f : Z \longrightarrow W$ be a holomorphic map between complex manifolds.
Assume that $(X_s)_{s \in S}$ is an analytic family of $n-$cycles in $Z$ and $Y$
a codimension $(n+1)-$cycle in $W$ such that

\smallskip

\item{1)} $\big(f_* (X_s)_{s \in S}\big) $ is an analytic family of $n-$cycles in 
$W$\footnote{(*)}{This is so iff $ \forall s \in S$, $ f|_{|X_{s}|}$ is proper and
generically finite on its images ; see [B.1] ch.IV.} 

\smallskip

\item{2)} $ f^*(Y)$ is a codimension $(n+1)$ cycle in $Z$\footnote{(**)}{This is
satisfied iff $f^{-1}(|Y|)$ has pure codim $n+1$ in $Z$ ; see [B.1] ch.VI \S  3.}
and satisfies $[AP]$ for the family $(X_s)_{s \in S}$.

\smallskip

\noi Then $Y$ satisfies $[AP]$ for the family $\big(f_*(X_s)\big)_{s \in S}$
and we have
$$\vbox{\boxit{5pt}{$\sum_{f^{*}(Y)}=\sum_Y$}}$$

\noi {\bf Proof.}

\noi Let $ F:={\rm id}_S \times f:S \times Z \longrightarrow S \times W.$ Condition
1) implies that $ F :|X| \longrightarrow |\widehat{X}|,$ where $\widehat{X}$ is the
graph of the family $\big(f_*(X_s)\big)_{s \in S},$ is proper and generically finite on
each fiber over $S.$ So Condition $[AP]$ for $f^*(Y)$ and the family $(X_s)_{s \in S}$
gives $[AP]$ for $Y$ and the family $\big(f_*(X_s)\big)_{s \in S}$ \big(recall
that $ F:|X|\longrightarrow|\widehat{X}|$ is surjective from 1)\big).

\smallskip

\noi The relation  $C^{W}_{\widehat{X}/S}=F_*(C^Z_{X/S})$ is true in relative holomorphic cohomology \big(see [B.1] ch.IV and [B.2]\big) and then still true in relative Deligne cohomology (see [K], the semi-purity properties); for this, see the  Appendix [C].

\noi Thanks to the projection formula (see Appendix 3), the section of the sheaf $ \underline{I\!\!H}^2_T \big(\Bbb Z(1)_{\cal D})$ associated to $ \Sigma_Y$ is given by :
$$\matrix{
Tr_{/S} \big(F_* (C^{S\times Z}_{X/S}) \cup C^{S\times W}_{S\times Y}\big) &=& Tr_{/S} \big(F_* (C^{S\times Z}_{X/S} \cup F^* C^{S\times W}_{S\times Y})\big) \cr
&=& \widehat{Tr}_{/S} (C^{S\times Z}_{X/S} \cup F^* C^{S\times W}_{S\times Y}) \hfill \cr}$$
because $Tr_{/S}  F_* = \widehat{Tr}_{/S}$ (see Appendix [C]). 

\noi So we conclude that we have the same section of $\underline{I\!\!H}^2_T \big({\Bbb Z}
(1)_{\cal D}\big)$ and $\Sigma_{f^{*}Y} = \Sigma_Y$. 

\smallskip

\noi {\bf Remark.}

\noi Using Property 2) to normalize $S$ and Property 3) to localize at generic points of $T=
|\Sigma_Y|$ we can reduce the proof of the formula $\Sigma_{f^{*}(Y)} = \Sigma_Y$ to the
equality of corresponding Weil divisors. This can be deduced from the projection formula
(then $X$ and $\widehat{X} = F_* X$ are seen as absolute cycles in $ S \times Z$
and $S \times W$, with $S$ smooth). This avoids to use compatibility for the relative
trace in relative Deligne cohomology with direct images (proved in the Appendix [C]).

\smallskip

\noi {\bf 5) Inverse image for cycles.}

\noi Consider  a holomorphic map $f:Z \longrightarrow W$ between complex
manifolds, $(X_s)_{s \in S}$  an analytic family of cycles in  $W$ such that the family $(f^* X_s)_{s \in S}$ is an analytic family of $n-$cycles
in $Z$ \big(see [B.1] ch.VI \S 3 for a definition\big) and $Y$ an  $(n+1)-$codimensional
cycle  in $Z$ such that $f_* Y$ is well defined in $W$.

\noi Assume that $Y$ satisfies $[AP]$ for the family $\big(f^*(X_s)\big)_{s \in S}$,
 then  we want to show  that $f_* Y$ satisfies $[AP]$ and that
$$\vbox{\boxit{5pt}{$\sum_Y=\sum_{f_*Y}$}}$$

\noi {\bf Proof.}

\noi Again Condition $[AP]$ for the cycle $f_* Y$ is easy to verify because
$f\big||Y|$ is proper and surjective on $|f_*(Y)|$. Then using as before the projection
formula (Appendix [B]) and the compatibility of the (relative) trace with direct image we 
have
$$\matrix{
Tr_{/S} \Big(F_*\big(F^*(C^{S\times W}_{X/S}) \cup C^{S\times Z}_{S\times Y}\big)\Big) &=& Tr_{/S} \big(C^{S\times W}_{X/S} \cup F_* C^{S\times Z}_{S\times Y}
\big) \hfill \cr
&=& \widehat{Tr}_{/S} \big(F^*(C^{S\times W}_{X/S}) \cup C^{S\times Z}_{S\times Y}\big). \cr}$$

\noi Thus the corresponding sections of the sheaf $\underline{I\!\!H}^2_T \big(
{\Bbb Z}(1)_{\cal D}\big) $ coincide. $\hfill \blacksquare$ 

\smallskip

\noi {\bf 6) Parameters on $Y$.}

\noi Now we consider an analytic family $(Y_v)_{v \in V}$ of $(p-1)-$cycles in $Z$
parametrized by a reduced complex space $V$ and such that, for any $ v \in V,$ the
cycle $Y_v$ satisfies $[AP]$ for the family $(X_s)_{s \in S}$ in $Z$.

\noi Then the family $(\Sigma_{Y_{v}})_{v \in V}$ is a flat family of Cartier divisors in
$S$ or an analytic family of cycles ; for Cartier divisors it is equivalent ( see [B.1] ch. V).

\noi The proof here is a repetition of the absolute case where $V$ is a point;  the absolute  fundamental class of $Y$ in $Z$ being replaced  by the relative fundamental class of the family
$(Y_v)_{v \in V}$.

\smallskip

\noi Then, using $S \times V$ as a parameter space, we conclude that there exists a
Cartier divisor $\Sigma_{\cal Y}$ in $S \times V$ such that (using base change) for any $v
\in V,  \Sigma_{\cal Y}|_{S \times \{v\}} = \Sigma_{{ Y}_{v}}$.

\smallskip

\noi so $\Sigma_{\cal Y}$ is $V-$flat in $S \times V$ and the proof is complete. $\hfill
\blacksquare$

\vskip 1cm

\noi Now to complete the proof of Theorem 1 we shall show that Conditions 1) to 6)
characterize our construction :

\smallskip

\noi The first remark is that, using 2) and 3), our construction is uniquely determined by
the case $S$ smooth and $|p_{X}^*Y|$ {\bf finite} on $|\Sigma_{Y}|.$ Then using 6) we can move
locally $Y$ in order to have $|Y|$ smooth around $|Y| \cap |X_{s_{0}}|$ (see
  Appendix II.A  for details on this argument). In this case we can assume that the cycle
$Y$ underlies a complete intersection ideal (locally in $Z)$ and using 4) we can
reduce the situation to $Z = {\Bbb C}^{n+1}$ and $Y=k\{0\}$\footnote{(*)}{For this
argument see [BMg1] p. 831.}, $(X_s)_{s \in S}$ being an analytic family of hypersurfaces
near $0$ in ${\Bbb C}^{n+1}$. Finally, using 5) to cut with a smooth curve through
$\{0\},$ we are reduced to the case $Z= {\Bbb C}$ and $Y=k\{0\}$ which is determined
by  1). \hfill $\blacksquare$\smallskip

\noi {\bf Remark 1.}

\noi As we have seen in the end of the proof above, the condition 5) allows to cut
the cycles $(X_s)_{s \in S}$ with a submanifold containing $Y$ (when the intersections
have the expected dimension).

\smallskip

\noi {\bf Remark 2.}

\noi Assume that $Y$ and $Y'$ are admissible poles which are, as $(n+1)-$codimensional
cycles, rationally equivalent (i.e  $\exists  {\cal Y} \subset {\Bbb P}_1 \times Z$ codimension
$n+1$ cycle in ${\Bbb P}_1 \times Z,$ equidimensional on ${\Bbb P}_1,$ such that
${\cal Y}_0=Y$ and ${\cal Y}_{\infty} =Y'$). Then for $S$ compact the corresponding Cartier
divisors $\Sigma_Y$ and $\Sigma_{Y'}$ in $S$ are linearily equivalent. This is an
easy consequence of [K] because the holomorphic map ${\Bbb P}_1 \longrightarrow H^1(S,
{\cal O}^*_S)$ obtained from the analytic family $({\cal Y}_t)_{t \in {\Bbb P}_{1}}$ has to
be constant.\smallskip\noindent
Another easy consequence of [K] is the result concerning the Abel's theorem (or problem b) ) obtained by  Griffiths in [Gr].
\vfill\eject
\noi {\tite 4. Intersection number of the incidence divisor with  a curve.}

\smallskip

\noi We consider the situation of the Theorem . Let ${\cal C}$ be a curve in $S$ 
and assume

\smallskip

\item{1)} $ |{\cal C}| \cap |\Sigma_{ Y}| $ is finite

\smallskip

\item{2)} for any $ \sigma \in |{\cal C}| \cap |\Sigma_Y|$, we have $|X_{\sigma}| \cap |Y| $ is finite.

\smallskip

\noi Then the following proposition gives a very simple formula to compute $ \# ({\cal C}. \Sigma_Y)$.
Note  that this intersection product is well defined because ${\cal C}$ inducesa homology class and  $\Sigma_Y$ a  cohomology class in $S$.

\smallskip

\noi {\bf Proposition 1.}

\noi \it{ Let\smallskip\noindent
{$\bullet$} $(X_s)_{s \in S} $ be an analytic family of $n-$cycles in a complex manifold
$Z$\par\noindent
{$\bullet$} $Y$ be a $n+1-$ codimensional cycle in $Z$ which is an admissible pole for
the family $(X_s)_{s \in S}$ (see the condition $[AP]$ § 2)\par\noindent
{$\bullet$} $ {\cal C} \build \hookrightarrow_{}^{j} S$ be a curve satisfying conditions 1) and
2) above\par\noindent
{$\bullet$} $X_{\cal C}$ be  the graph of the family obtained by pull back by $j$
from $(X_s)_{s \in S}$\par\noindent
{$\bullet$} $p : S\times Z\rightarrow Z$ the canonical projection.\smallskip\noindent
Then there exist a neighbourghood $W$ of $\; \bigcup_{\sigma \in |{\cal C}| \cap |\Sigma_Y|}X_{\sigma} \cap Y\; $ for which  the direct image $p'_* (X_{\cal C})$ (where $p'$ is the restriction of $p$ to $S\times W$) of the cycle $X_{\cal C}$
is well defined  and we have :
$$\# \big(p_* (X_{\cal C}). Y) = \# ({\cal C}.\Sigma_Y)$$}
\rm\smallskip
\noi Remark that in $W$ the intersection $p'_* (X_{\cal C}).Y$
is well defined because these cycles have dimension and codimension $n+1$ in $W$
and finite intersection by Assumptions 1) and 2).

\smallskip

\noi {\bf Proof.}
\noi The assumptions give us such a $W$. We denote $p$ the projection $p'$.
\noi Let $ \nu : {\cal C}' \longrightarrow {\cal C}$ be the normalization of ${\cal C}.$
The map $j \circ \nu : {\cal C}' \longrightarrow S$ is then an admissible pull back for
the pole $Y$ and the family $(X_s)_{s \in S}.$ So we are reduced to prove our 
formula when $S={\cal C}'$ is a smooth curve and all $X_{\sigma} \cap Y$ are finite.
But in this case, as $S$ is smooth, we have 
$$\Sigma_Y = \pi_* (X \build \cdot_{p}^{} Y)$$ 
\noindent where $\pi$ is the projection of $S\times Z$ on $S$.
\noi Using the projection formula and the fact that, for finite cycles, the intersection number  is preserved by direct image, we deduce that  $\#({\cal C}.\Sigma_Y)=\#(X_{\cal C}. p^* Y)=\#(p_* X_{\cal C}.Y)$.
$\hfill \blacksquare$

\vfill\eject
\noi {\bf Example.} (very classical)

\noi Let $Y$ be a codimension $n+1$ cycle in ${\Bbb P}_N$ and let $(X_s)_{s \in S}
$ be the universal family of $n-$planes in ${\Bbb P}_N$ \big(so $S= {\rm Gr}(n+1,
N+1)$ is the Grassmann manifold of $n-$plane in $ {\Bbb P}_N\big).$ Let ${\cal C}$
be the line of $n$ planes in a $(n+1)-$planes $\Pi$ containing a given $(n-1)-$plane
$P$  in $\Pi$.

\noi Assume $P \cap Y = \emptyset$ and $Y\cap \Pi$ finite. Then $\#(\Sigma_Y \cap
{\cal C})=\#(Y\cap \Pi)={\rm deg} Y$.

\noi So the intersection number of the Chow divisor $\Sigma_Y$ of $Y$ in Gr$(n+1,
N+1)$ with such a generic "line" is the degree of $Y$ in ${\Bbb P}_N$.

\smallskip

\noi {\bf{Remark.}}\rm\par\noindent
 If we allow that, for some $\sigma$ in $|{\cal C}| \cap |\Sigma_Y|, 
|X_{\sigma}| \cap |Y|$ is not finite (but necessarily compact) the right handside of the
previous formula is not defined. One can try to replace it by ${\rm Tr}_W (C^W_{p_{*}(X_{\cal C})}
\cup C^W_Y)$ which is well defined and give the same number in the finite intersection
case.

\noi It is possible to extend the proposition with this formulation as soon as it is possible
to move the curve ${\cal C}$ in $S$ in an analytic family $({\cal C}_t)_{t \in T}$
such that ${\cal C}_{t_{0}} = {\cal C}$ and that for generic $t$ we have $|X_{\sigma}|\cap |Y|$ is finite  for all
$\sigma \in |{\cal C}_t| \cap |\Sigma_Y|$ \footnote{(*)}{for instance
if $|{\cal C}| \cap |\Sigma_Y|$ is in Reg $(S)$ this will be possible.} : then, using the stability of analytic family by direct image and the existence of relative fundamental class for analytic family of cycles in complex manifold, it is enough 
to prove, after weak normalization of $T$, that the family $(X_{{\cal C}_{t}})_{t \in T}$ is analytic. 

\smallskip

\noi In this situation we shall set 
$$\matrix{
&\# ({\cal C}_{t_{0}}.\Sigma_Y) &=& \# ({\cal C}_{t}.\Sigma_Y) \hfill& \hbox{for 
generic} \quad t \quad \hbox{(near}\quad t_0) \hfill \cr
\noalign{\vskip 0,3cm}
{\rm and} & \#(X_{{\cal C}_{t}}. Y) \hfill &=& {\rm Tr}_W (C^W_{p_{*} X_{{\cal C}_{t}}}
\cup C^W_Y) \hfill& \cr
\noalign{\vskip 0,3cm}
&&=&{\rm Tr}_W (C^W_{p_{*} X_{{\cal C}_{t_{0}}}} \cup C^W_Y) \cr}$$

\noi because $t \longrightarrow {\rm Tr}_W(C^W_{p_{*}X_{{\cal C}_{t}}} \cup C^W_Y)$ will
be weakly holomorphic \big(existence of a relative fundamental class for the analytic family
$p_*(X_{{\cal C}_{t}})_{t \in \widehat{T}}$ in $W$ where $ \widehat{T}$ is the
weak normalization of $T\big)$ and with value in $\Bbb N$ for generic $t \in T$.

\noi To complete this extention of the Proposition it is enough  to prove the following

\noi {\bf Proposition 2.}

\noi \it{ Let $ \widetilde{Z} \build \longrightarrow_{}^{\pi} S$ be a geometrically flat map
relative to the geometric weight $\widetilde{X} \subset S \build \times_{S}^{}
\widetilde{Z}$ \big(see [BMg1] p.12 \big). Let $({\cal C}_{t})_{t \in T}$ be an 
analytic family of cycles in $S.$ Define for each $t \in T,  \widetilde{X}_{{\cal C}_{t}}$
as the cycle in $\widetilde{Z} \simeq \{t\} \times \widetilde{Z}$ defined as the graph
of the family of cycles in $\widetilde{Z}$ given by the graph of the base change $
{\cal C}_t \hookrightarrow S$ (see below for a more precise definition).

\noi Let $ \widehat{T} \longrightarrow T$ be the weak normalization of $T.$ Then
$(\widetilde{X}_{{\cal C}_{t}})_{t \in \widehat{T}}$ is an analytic family of cycles
in $\widetilde{Z}$.}\rm

\smallskip

\noi {\bf Proof.}

\noi Let us, first, be more precise on the definition of the cycles $\widetilde{X}_{{\cal C}_{t}}$
: let $Y=\dis \sum_{\alpha \in A} m_{\alpha} Y_{\alpha}$ be a cycle in $S$ with
$Y_{\alpha}$ irreducible and locally finite, $m_{\alpha} \in \Bbb N^*.$ We define
$\pi^*(Y)$ in $\widetilde{Z}$ as follows : let ${\cal Y}_{\alpha}$ be the cycle in
$\widetilde{Z}$ which is the graph of the analytic family defined by the base change
$Y_{\alpha} \hookrightarrow S$ from the family $(\widetilde{X}_s)_{s \in S}$ of
fibers of $\pi : (\widetilde{Z},\widetilde{X}) \longrightarrow S.$

\noi Then let $\pi^*(Y):=\dis \sum m_{\alpha} {\cal Y}_{\alpha}$.

\noi So in the Proposition $\widetilde{X}_{{\cal C}_{t}}:=\pi^*({\cal C}_t)$.

\smallskip

\noi Now, using the fact that the change of projection (in a scale) is always meromorphic
and continuous \big([B.1] Th. 2 p.42 ; we can also use proposition 3 of ch. III \S 3\big) it
is enough to prove that for any $t_0 \in T$ and any $z_0 \in \pi^{-1}\big(|{\cal C}_{t_{_{0}}}|\big)$
we can find a scale in a neighbourghood of $z_0$ adapted to $\pi^*({\cal C}_{t_{0}})$
such that the family of branched  covering associated to $\pi^*({\cal C}_t)_{t \in T}$ is
analytic near $t_0$.

\smallskip

\noi For that purpose, choose a scale adapted to $ {\cal C}_{t_{0}}$ near $\pi(z_0)$
given by a local embedding
$$S \hookrightarrow V \times B'$$

\noi and a scale adapted to $\widetilde{X}_{\pi(z_{0})}$ near $z_0$ given by a local
embedding :
$$\widetilde{Z} \hookrightarrow S \times U \times B \hookrightarrow V \times U \times B
\times B'.$$

\noi So we have holomorphic maps :
$$\matrix{
& f & :& S' \times U & \longrightarrow &{\rm Sym}^l(B) \cr
&g &:& T' \times V & \longrightarrow &{\rm Sym}^k(B') \cr} $$

\noi defining respectively the fibers of $\pi$ near $\pi(z_0) \in S'$ and the
family $({\cal C}_t)_{t \in T}$ near $t_0 \in T'$ (here $S'$ and $T'$ are
"small" open sets in $S$ and $T$).

\noi Remark that, by our assumptions, $f$ and $g$ are isotropic \big(see [B.1] ch.
2\big).

\noi Then the family $\pi^*({\cal C}_t)_{t \in T'}$ is defined in the scale given by
$\widetilde{Z} \hookrightarrow (V \times U) \times (B \times B') $ (near $z_0)$ by
the map 
$$T' \times V \times U \longrightarrow {\rm Sym}^{kl} (B\times B')$$

\noi obtained as follows : from $g$ we have an holomorphic map $G:T'\times V\longrightarrow
{\rm Sym}^k(V\times B')$ where $G(t,v)=\big((v,b'_1),\cdots,(v,b'_k)\big)$,  if $g
(t,v)$ is the $k-$uple $(b'_1,\cdots,b'_k).$ But, by assumption, $(v,b'_{j}) \in S$ for
$j \in [1,k]$. So $G$ factors through the holomorphic map
$$G_0 : T' \times V \longrightarrow {\rm Sym}^k(S').$$

\noi Now, using $G_0$ and $f$ we obtain an analytic  map
$$T' \times V \times U \longrightarrow {\rm Sym}^k(S' \times U) \build \longrightarrow_{}^{{\rm
Sym^{k} F}} {\rm Sym}^{kl} (B \times B') $$

\noi where $F:S' \times U \longrightarrow {\rm Sym}^l (B \times B')$ is given by
$$F (s,u) = \big((b_1,p(s)\big),\cdots,\big(b_l,p(s)\big)$$

\noi if $f(s,u)=(b_1,\cdots,b_l)$ and $p:S' \longrightarrow B'$ is the composition of the
embedding $S' \hookrightarrow V \times B'$ and the projection on $B'$.

\noi So the map $T' \times V \times U \longrightarrow {\rm Sym}^{kl} (B \times B')$ is
holomorphic and this proves Proposition 2. $\hfill \blacksquare$

\bigskip

\noi {\bf Remark.}

\noi It is not clear (an may be not true ) that the map $T' \times V \times U \longrightarrow
{\rm Sym}^{kl} (B \times B')$ is $T'-$isotropic when $f$ and $g$ are isotropic.
That is the reason why we have to normalize weakly $T$ in our conclusion.

\noi Of course, because the weak normalization is an homeomorphism, this is irrelevant in
our generalisation of the proposition 1.

\noi To conclude, let us state the simple case where this generalization works :

\bigskip

\noi {\bf Corollary.}

\noi\it{ In the situation of Proposition 1, assume  

\noi 2') for any $\sigma \in |{\cal C}| \cap
|\Sigma_Y|$, either  $ |X_{\sigma}| \cap|Y|$ is finite or $\sigma$ is a smooth point in
$S$.

\noi Then Hypothesis  1) and 2') imply 
$$\# ({\cal C}.\Sigma_Y)= {\rm Trace}_W (C^W_{p_{*}(X_{\cal C})} \cup C^W_Y)$$}\rm

\centerline{---------------------------------------}
\vfill \eject

\centerline{\titre Appendix}

\vskip 1cm

\noi {\tite I.} 

\smallskip

\noi The formulas we want to prove in this Appendix (intersection, projection and compatibility
between direct image and trace in Deligne cohomology) are consequences of

\smallskip

\item{1)} the existence in Deligne cohomology of a push-foward and a pull-back for an holomorphic map between 
complex manifolds.
\smallskip

\item{2)} the corresponding formulas in holomorphic cohomology

\smallskip

\noi and we conclude using the injectivity result of [K] corollary 2 p.295.

\noi So our first purpose is to prove 1).

\noi Let $f:Z\longrightarrow W$ be a holomorphic map between connected complex manifolds
and let $d={\rm dim}_{\Bbb C} W-{\rm dim}_{\Bbb C} Z\quad (d \in \Bbb Z)$.

\noi To build a push-forward for  $k \in {\Bbb N}$
$$f_*^{\cal D} : Rf_{\rm !}  {\Bbb Z}(k)_{{\cal D}, Z} \longrightarrow {\Bbb Z}(k+d)_{{\cal D}, W}
[2d] \leqno (1)$$ 

\noi we shall use the holomorphic and topological push-forwards
$$\eqalign{
& f^h_* : Rf_{\rm !}  \Omega^{\bullet}_{Z} \longrightarrow \Omega^{{\bullet}+d}_W [d] \cr
& f^t_* : Rf_{\rm !} {\Bbb Z}_{Z}  \longrightarrow {\Bbb Z}_{W} [2d] \cr} \leqno (2)$$

\noi which are defined as follows :

\smallskip

\noi Using the factorization of $f$ by its graph, it is enough to handle the case
where $f$ is the embedding of a closed submanifold and when $f$ is the projection
of a product (so $Z=X \times W$ with ${\rm dim}_{\Bbb C} X=d)$.

\smallskip

\noi In the first case $d={\rm codim}_Z W,$ we have
$$\eqalign{
& f^h_* : \Omega^{\bullet}_{Z} \longrightarrow \underline{H}^d_Z (\Omega^{{\bullet}+d}_W) \cr
{\rm and} \qquad
& f^t_* : {\Bbb Z}_{Z} \longrightarrow \underline{H}^{2d}_Z ({\Bbb Z}_W) \cr }$$

\noi which are quasi isomorphisms.

\noi In the case of the projection $ Z=X \times W \longrightarrow W$ we have the usual
"integration" maps.
$$\eqalign{
& f^h_* : R^d f_{\rm !}  \Omega^{{\bullet}+d}_Z \longrightarrow \Omega^{\bullet}_W \cr
& f^t_* : R^{2d} f_{\rm !} { \Bbb Z}_{Z} \longrightarrow {\Bbb Z}_{W} \cr }$$

\noi The push-forward (1) is now built by a decreasing induction on $k$.

\noi For $k$ large enough we have, by the exactness of the holomorphic de Rham complex,
a quasi-isomorphism
$${\Bbb Z}(k)_{\cal D} \simeq ({\Bbb C}/{\Bbb Z}) [-1]$$

\noi on any complex manifold. So the topological push-forward $f^t_*$ with coefficient
$\Bbb C/\Bbb Z$ is enough do define (1) for large $k$.

\noi To define (1) for $k$ assuming that it is already built for $k+1$, we  consider
the exact sequence of complexes on $Z$ and $W$
$$0 \longrightarrow \Omega^k [-k-1] \longrightarrow {\Bbb Z}(k+1)_{\cal D} \longrightarrow
{\Bbb Z}(k)_{\cal D} \longrightarrow 0$$

\noi and the diagram

\def\hfl#1#2{\smash{\mathop{\hbox to 6mm{\rightarrowfill}}\limits^{\scriptstyle#1}_{
\scriptstyle#2}}}
\noi $\diagram{
Rf{\rm !}  \Omega^k_Z [-k-1] & \hfl{}{}& Rf{\rm !}  {\Bbb Z}(k+1)_{{\cal D}, Z} &
\hfl{}{} & Rf{\rm !}  {\Bbb Z}{k}_{{\cal D}, Z} & \hfl{}{} \cr
\vfl{}{f^{h}_{*}} && \vfl{}{f^{\cal D}_*}&& \vfl{}{} \cr
\Omega^{k+d}_W (d-k-1][2d] & \hfl{}{} & {\Bbb Z}(k+d+1)_{{\cal D}, W}[2d]& \hfl{}{}&
{\Bbb Z} (k+d)_{{\cal D},W}[2d] & \hfl{}{} \cr } $

\noi from which we deduce the definition of $f^{\cal D}_*$ for $k$ (we assume inductively
on $k$ the commutavity of this diagram).\smallskip\noindent
We can claim that\smallskip\noindent
1) The functor $Rf_{!}$ transforms distinguished  triangles into distinguished triangles.\par\noindent
2) If $(A, B, C)$ and $(A', B', C')$ are two distinguished triangles with two morphisms $u : A \rightarrow A'$ and $u' : B\rightarrow B'$ such that the square $(A, B, A', B')$ is commutative then we have a morphism between these triangles. 

\noi The construction for the pull-back $\qquad (\forall k \in \Bbb N)$
$$f^*_{\cal D} : f^{\rm !}  {\Bbb Z}(k)_{{\cal D}, W} \longrightarrow {\Bbb Z}(k-d)_{{\cal D},
Z} [-2d] \leqno (3) $$

\noi works in the same way :

\noi We use the holomorphic and topological pull-back
$$ \eqalign{
& f^*_h : Rf^{\rm !}  \Omega^{\bullet}_W \longrightarrow \Omega^{{\bullet}-d}_W [-d] \cr
& f^*_t : f^{\rm !}  {\Bbb Z}_{W} \longrightarrow {\Bbb Z}_{Z} [-2d]} \leqno (4) $$

\noi defined, as before, by factorizing $f$ through its graph. Then for $f$ the 
inclusion of $Z$ as a codimension $d$ closed submanifold of $W$ we have
$$\eqalign{
&f^*_h :  \underline{H}^{d}_Z (\Omega^{{\bullet}+d}_W) \longrightarrow \Omega^{\bullet}_Z \cr
{\rm and} \qquad
& f^*_t : \underline{H}^{2d}_Z  ({\Bbb Z}_W) \longrightarrow {\Bbb Z}_Z \cr} $$

\noi which are the holomorphic and topological residue maps (they are inverse of $f^h_*$
and $f^t_*$ respectively ).

\noi In the projection case we have
$$\eqalign{
& f^*_h : f^{\rm !} \Omega^{\bullet}_W = f^* \Omega^{\bullet}_{W} \otimes \Omega^d_{{W\times Z}/W} [d] \longrightarrow
\Omega^{{\bullet}+d}_Z [d] \cr
{\rm and} \qquad
& f^*_t : f^{\rm !} {\Bbb Z}_{W}=f^* {\Bbb Z}_{W} \otimes {\Bbb Z}_{Z\times W} [2d] \longrightarrow {\Bbb Z}_{Z}
[2d]. \cr }$$

\noi Then we proceed, as before, by a decreasing induction on $k$ to construct the
pull-back (3).

\bigskip

\noi {\bf Remark.}

\noi By construction, it is obvious that (1) and (3) are compatible with (2) and (4).

\vskip 1cm

\noi {\tite II.}

\smallskip

\noi Our aim now is to describe more precisely the nice behaviour of fundamental classes in Deligne cohomology
(absolute and relative case) for some simple operations. These results are not "new" but
we were unable to give references for these results in our local analytic context.

\bigskip

\noi {\bfdouze A. Intersection.}

\noi We shall prove the following

\smallskip

\noi {\bf Proposition 1.}

\noi {\it Let \par\noindent
{$\bullet$} $Z$ be a complex manifold\par\noindent
{$\bullet$} $S$ and $T$ complex analytic  reduced spaces\par\noindent
{$\bullet$} $(X_s)_{s \in S}$ (resp. $(Y_t)_{t \in T}$)
  analytic family of pure $p$ (resp. $q$)- codimensional cycles in $Z$ parametrized by the reduced analytic
space $S$ (resp.  $T$).\par\noindent
$\bullet$  $C^Z_{X/S}$ and $C^Z_{Y/T}$ the
relative fundamental classes in Deligne cohomology for these families\par\noindent
{$\bullet$} $p_{1} : S\times T\times Z \rightarrow S\times Z$ (resp.$q_{1} : S\times T\times Z \rightarrow T\times Z$) the canonical projection. \par\noindent

\noi Assume that for any $(s,t) \in S \times T$ the closed analytic set $|X_s| \cap
|Y_t|$ has the expected codimension.\smallskip\noindent
 Then,  $(X_{s}.Y_{t})_{(s,t)\in S\times T}$  is an analytic family of cycles in $Z$ and its  relative fundamental class in Deligne cohomology is given by the cup product $p^{*}_{1}C^Z_{X/S} \cup
q^{*}_{1}C^Z_{Y/T} $ in $ \underline{I\!\!H}^{2(p+q)}_{p_{1}^{-1}(|X|) \cap q_{1}^{-1}(|Y|)} ({\Bbb Z}(p+q)_{{\cal D}/S
\times T})$ }
\vfill\eject
\rm
\noi {\bf Proof.}

\noi We shall proceed by successive reductions

\item{} {\bf First reduction}\par\noindent
First, let $X_{1}:= p_{1}^*(X)$ and $Y_{1}:= q_{1}^*(Y)$.\par\noindent
Theorem. 10
(local) of [B.1] ch.VI says us that the family
$(X_{s}.Y_{t})_{(s,t)\in S\times T}$ is analytic.
\noi For our purpose, it suffises  to work in holomorphic cohomology because we have injective morphisms of
sheaves :
$$\matrix{
\hfill \underline{I\!\!H}^{2p}_{|X|} \big({\Bbb Z}(p)_{{\cal D}/S}\big) & \longrightarrow &
\underline{H}^{p}_{|X|} (\Omega^p_{S \times Z/S}) \hfill \cr
\hfill \underline{I\!\!H}^{2q}_{|Y|} \big({\Bbb Z}(q)_{{\cal D} /S}\big) & \longrightarrow &
\underline{H}^{q}_{|Y|} (\Omega^q_{T \times Z/T}\big) \hfill \cr
\underline{I\!\!H}^{2(p+q)}_{|X_{1}| \cap |Y_{1}|} \big({\Bbb Z}(p+q)_{{\cal D} /S \times T}\big) &
\longrightarrow & \underline{H}^{p+q}_{|X_{1}| \cap |Y_{1}|} (\Omega^{p+q}_{S \times T \times Z/S
\times T} )\cr}$$
\noindent \big( see [K] p. 320\big)

\item{} {\bf Second reduction}

\noi We can assume $S$ and $T$ are points, because the sheaf $ \underline{H}^{p+q}_{|X_{1}| \cap
|Y_{1}|} (\Omega^{p+q}_{S \times T \times Z/S \times T})$ has no torsion over ${\cal O}_{S \times T}$

\smallskip

\item{} {\bf Third reduction}

\noi Using diagonal trick, we can assume that $Y$ is a smooth (closed) submanifold in $Z$.

\smallskip

\noi Now the case where $X$ and $Y$ are
smooth and transversal along $ |X| \cap |Y|,$ is trivial. We denote this intersection $X\bullet Y$.

\smallskip

\noi We are reduced to consider the case where $X=|X|$ is irreducible of dimension $n$
and where $Y$ is smooth of codimension $n-q$ in $Z$ an open ball in $ \Bbb C^N.$
Of course we assume that $X \cap Y$ is of pure dimension $q.$ Then we consider the
family of $Y_s = Y+s$ (translated of $Y$) where $s \in \Bbb C^N$ is small
enough, say in $S$. Then for any given open set \footnote{(*)}{see [B.1] proposition 2 p. 138 for a more precise proof
of this.} in $Z$, there exists an open dense set in $S$ such that the intersection  $X.Y_{s}$ is smooth and transversal. So $C^Z_X \cup C^Z_{Y/S}$ induces, for an open dense set in
$S,$ the fundamental class of $X \cap Y_s.$ So $C^Z_{X/S} \cup C^Z_{Y/S}$ has to
be the relative fundamental class of the family $ (X \cap Y_s)_{s \in S},$ which is
analytic by [B.1] Th. 10 (local) and so has a relative fundamental class by [B.2]. So we
conclude that for $s=0$ we have 
$$C^Z_X \cup C^Z_Y = C^Z_{X \cap Y}. $$
$\hfill \blacksquare$\smallskip
\noi {\bfdouze B. Projection formula.}

\smallskip

\noi We want now to prove the following
\vfill\eject
\noi {\bf Proposition 2.}

\noi {\it Let $f:Z \longrightarrow W$ be a holomorphic map between complex manifolds of
dimension $d$ and $\delta$ respectively. Let $X$ be a $n-$cycle in $Z$
such that $f|_{|X|}$ is proper and generically finite on its image ; let $Y$ be a
$(n+q)-$codimensional cycle in $W.$ Assume that $X \build {\cdot}_{f}^{} Y$ is
well defined\footnote{(**)}{\rm see [B.1] chap.6 \S 3. This hypothesis means that in $Z \times
W,  Z \times Y$ meets the

\noi graph of $f|_{|X|}$ in the expected dimension.}. Then we can define $f_* X$ in
$W$ and also $(f_* X) . Y$. The direct image $f_*(X \build \cdot_{f}^{} Y)$
is also well defined and the projection formula says 
$$f_* (X \cdot_{f}^{} Y)=(f_* X) \cdot Y. \leqno (P_0) $$

\noi The corresponding formula in Deligne cohomology is the equality
$$f_* \big(C^Z_X \cup f^*(C^W_Y)\big) = f_*(C^Z_X) \cup C^W_Y \leqno (P)$$

\noi where the map $f^*: \underline{I\!\!H}^{2(n+q)}_{Y} \big(W,  {\Bbb Z}(n+q)_{\cal D}
\big) \longrightarrow {I\!\!H}^{2(n+q)}_{f^{-1}(Y)} (Z, {\Bbb Z}(n+q)_{\cal D}\big)$ is the
pull back and the maps
$$\matrix{
f_*: \underline{I\!\!H}^{2(d-q)}_{|X| \cap f^{-1}(|Y|)} \big(Z, {\Bbb Z}(d-q)_{\cal D}\big) &
\longrightarrow & \underline{I\!\!H}^{2(\delta -q)}_{f|X| \cap |Y|} \big(W, \Bbb Z(\delta
-q)_{\cal D}\big) \cr
\hfill f_*: \underline{I\!\!H}^{2(d-n)}_{|X|} \big(Z, {\Bbb Z}(d-n)_{\cal D}\big) &
\longrightarrow & \underline{I\!\!H}^{2(\delta -n)}_{f|X|} \big(W, 
{\Bbb Z}(\delta -n)_{\cal D}\big) \hfill \cr}$$

\noi are the trace maps (described below via the holomorphic cohomology).}

\bigskip

\noi {\bf Proof.}

\noi As before it is enough to prove Formula $(P_0)$ in holomorphic cohomology. Recall
that the trace map, for a closed set $F$ such that $f|_F$ is proper,
$$f_* : H^{d-\alpha}_F (Z, \Omega^{d-\alpha}_Z) \longrightarrow H^{\delta - \alpha}_{f(F)}
(W, \Omega^{\delta - \alpha}_W)$$

\noi is given by direct image of hyperfunctions of type $(d-\alpha,  d-\alpha),
\overline{\partial}-$closed with support in $F$.

\noi Because the part $A$ and the previous remark we have only to prove that the formula 
$$C^W_{f_{*}(X)} = f_* (C^Z_X)$$

\noi is valid in holomorphic cohomology. We can use for that, the fact that the integration
current $[X]$ on $X$ represents the class $C^Z_X $ in $H^{d-n}_{|X|} (Z,
\Omega^{d-n}_Z).$ Now the formula
$$f_* [X] = [f_* X]$$

\noi between currents in $W$ is essentially a trivial change of variable formula. So the
proof is complete. $\hfill \blacksquare$
\bigskip

\noi {\bf Corollary.}

\noi {\it The Proposition extends in a straightforward manner in  the relative case $\big(
(X_s)_{s \in S} $ in $Z$ and $(Y_t)_{t \in T}$ in $W$ are now analytic families
of cycles\big).}

\vskip 1cm

\noi {\bfdouze C. Compatibility  between direct image and trace}

\noi We shall prove the following

\smallskip

\noi {\bf Proposition 3.}

\noi {\it Let $f : Z \longrightarrow W$ be a holomorphic map between complex manifolds of
dimension $d$ and $\delta$ respectively. Then we have the following commutative
triangle }


$$\xymatrix{{I\!\!H}^{2d}_{c}\big({\Bbb Z}(2d)_{{\cal D},Z}\big)\ar[rr]^{f_{*}}\ar[rd]^{\widetilde{Tr}}&&\ar[ld]_{Tr}{I\!\!H}^{2\delta}_{c}\big({\Bbb Z}(2\delta)_{{\cal D},W}\big) \\ 
&\Bbb{C}&}$$
{\it \noi where $ \widetilde{Tr}$ and $Tr$ are the trace maps.}

\bigskip

\noi {\bf Proof.}

\noi In fact if $H^{d}_{c} (Z, \Omega^{d}_{Z}) \build \longrightarrow_{}^{\int_Z} {\Bbb C} $ 
\big(resp. $H^{\delta}_c (W, \Omega^{\delta}_{W}) \build \longrightarrow_{}^{\int_W}{\Bbb C}\big)$ 
is the usual integration map, $\widetilde{Tr} $ (resp. $Tr$) is the composition of 
$\int_Z$ (resp. $ \int_W)$ with the connector of the long exact hypercohomology 
sequence of the short exact sequence of complexes
$$0 \longrightarrow \Omega^d [-d-1] \longrightarrow {\Bbb Z}(d+1)_{{\cal D}} \longrightarrow
{\Bbb Z}(d)_{{\cal D}} \longrightarrow 0.$$

\noi The commutativity of the triangle above is then a consequence of the
compatibility of the integration maps with the direct image : if $\varphi$ is a $
C^{\infty} (d,d)$ form with compact support on $Z, f_* \varphi$ is a $(\delta,
\delta)$ current with compact support in $W$ and we have 
$$\int_W  f_* \varphi := \langle f_* \varphi, 1 \rangle = \langle \varphi, f^*(1) \rangle
= \int_Z \varphi .$$
$\hfill \blacksquare$

\vskip 1cm

\noi {\bf Corollary.}

\noi {\it Proposition $3$ holds  in the relative setting.}
\bigskip\bigskip\bigskip

\vfill\eject {\centerline{\bfdouze  References}}

\bigskip

\item{[B1]} D. Barlet : {\it Espace analytique r\'eduit des cycles analytiques complexes
compacts d'un espace analytique r\'eduit.} Sem. F. Norguet, Lectures Notes in Mathematics,
{\bf 482}, (1975), p. 1-158.

\smallskip

\item{[B2]} D. Barlet : {\it Familles analytiques de cycles et classes fondamentales relatives.}
Sem. F. Norguet, Lectures notes in Mathematics, {\bf 807}, (1980), p. 1-24.

\smallskip

\item{[B3]} D. Barlet : {\it Faisceau $\omega^{\cdot}$ sur un espace analytique de dimension
pure.} Sem. F. Norguet, Lectures notes in Mathematics, {\bf 670}, (1978), p. 187-204. 

\smallskip

\item{[BMg1]} D. Barlet, J. Magnusson : {\it Int\'egration de classes de cohomologie m\'eromorphes
et diviseur d'incidence.} Ann. Sc. de l'ENS, s\'erie 4, {\bf t. 31}, fasc. 6, (1998), p. 811-
842. 

\smallskip

\item{[BMg2]} D. Barlet, J. Magnusson : {\it Transfert de l'amplitude du fibr\'e normal au
diviseur d'incidence.} J. Reine Angew. Math., {\bf 513}, (1999), p. 71-95.

\smallskip

\item{[D]} P. Deligne : {\it Le d\'eterminant de la cohomologie.} Contemporary Mathematics,
{\bf 67}, (1985), p. 95-177.

\smallskip

\item{[E]} R. Elkik : {\it Fibr\'es d'intersection et int\'egrales de classes de Chern.} Ann.
Sc. de l'ENS, s\'erie 4, {\bf t. 22}, (1989), p. 155-226.

\smallskip

\item{[EV]} H. Esnault, E. Viehweg : {\it Deligne-Beilinson cohomology} 
Perspectives. Math.
Academic Press, {\bf 4}, (1987), p. 43-92.

\smallskip

\item{[F]} W. Fulton : {\it Intersection theory.} Springer Verlag, Berlin, 1984.
\smallskip

\item{[Gr]} P. Griffiths: {\it Some results on algebraic cycles in algebraic manifolds.} 
Internat. Colloq., Tata Inst. Fund. Res., Bombay, (1968) pp. 93--191 Oxford Univ. Press.

\item{[K]} M. Kaddar : {\it Classe fondamentale relative en cohomologie de Deligne et
application.} Mathematische Annalen, {\bf 306}, (1996), p. 285-322.

\end\end